\documentstyle[12pt,amssymb]{amsart}
\def\psfig#1{}

\begin{document}

\newtheorem{thm}{Theorem}[section]
\newtheorem{cor}[thm]{Corollary}
\newtheorem{lem}[thm]{Lemma}
\newtheorem{prop}[thm]{Proposition}
\newtheorem{schw}{Schwarz Lemma} 
\newtheorem{sectl}{Sector Lemma}
\theoremstyle{remark}
\newtheorem{rem}{Remark}[section]
\newtheorem{notation}{Notation}

\numberwithin{equation}{section}
\newcommand{\thmref}[1]{Theorem~\ref{#1}}
\newcommand{\propref}[1]{Proposition~\ref{#1}}
\newcommand{\secref}[1]{\S\ref{#1}}
\newcommand{\lemref}[1]{Lemma~\ref{#1}}
\newcommand{\corref}[1]{Corollary~\ref{#1}} 
\newcommand{\figref}[1]{Fig.~\ref{#1}}

\renewcommand{\theschw}{}
\renewcommand{\thesectl}{}

\theoremstyle{definition}
\newtheorem{defn}{Definition}[section]

\newcommand{\QED}{\rlap{$\sqcup$}$\sqcap$\smallskip}

\newcommand{\di}{\partial}
\newcommand{\dibar}{\bar\partial}
\newcommand{\ra}{\rightarrow}
\newcommand{\hra}{\hookrightarrow}
\def\lra{\longrightarrow}

\def\ssk{\smallskip}
\def\msk{\medskip}
\def\bsk{\bigskip}
\def\noi{\noindent}
\def\nin{\noindent}
\def\lqq{\lq\lq}
\def\sm{\smallsetminus}
\def\bolshe{\succ}
\def\ssm{\smallsetminus}
\def\tr{{\text{tr}}}

\newcommand{\diam}{\operatorname{diam}}
\newcommand{\dist}{\operatorname{dist}}
\newcommand{\distM}{\operatorname{dist}_M}
\newcommand{\cl}{\operatorname{cl}}
\newcommand{\inter}{\operatorname{int}}
\renewcommand{\mod}{\operatorname{mod}}
\newcommand{\tl}{\tilde}

\newcommand{\orb}{\operatorname{orb}}
\newcommand{\HD}{\operatorname{HD}}
\newcommand{\supp}{\operatorname{supp}}
\newcommand{\id}{\operatorname{id}}
\newcommand{\length}{\operatorname{length}}
\newcommand{\width}{\operatorname{width}}
\newcommand{\dens}{\operatorname{dens}}
\newcommand{\meas}{\operatorname{meas}}

\newcommand{\Dil}{\operatorname{Dil}}
\newcommand{\Ker}{\operatorname{Ker}}
\newcommand{\tg}{\operatorname{tg}}
\newcommand{\codim}{\operatorname{codim}}
\newcommand{\isom}{\approx}
\newcommand{\esssup}{\operatorname{ess-sup}}

\newcommand{\SLa}{\underset{\La}{\Subset}}

\newcommand{\const}{\mathrm{const}}
\def\loc{{\mathrm{loc}}}

\newcommand{\eps}{{\epsilon}}
\newcommand{\De}{{\Delta}}
\newcommand{\de}{{\delta}}
\newcommand{\la}{{\lambda}}
\newcommand{\La}{{\Lambda}}
\newcommand{\si}{{\sigma}}
\newcommand{\Om}{{\Omega}}
\newcommand{\om}{{\omega}}

\newcommand{\AAA}{{\cal A}}
\newcommand{\BB}{{\cal B}}
\newcommand{\CC}{{\cal C}}
\newcommand{\DD}{{\cal D}}
\newcommand{\EE}{{\cal E}}
\newcommand{\EEE}{{\cal O}}
\newcommand{\II}{{\cal I}}
\newcommand{\FF}{{\cal F}}
\newcommand{\GG}{{\cal G}}
\newcommand{\JJ}{{\cal J}}
\newcommand{\HH}{{\cal H}}
\newcommand{\KK}{{\cal K}}
\newcommand{\LL}{{\cal L}}
\newcommand{\MM}{{\cal M}}
\newcommand{\NN}{{\cal N}}
\newcommand{\OO}{{\cal O}}
\newcommand{\PP}{{\cal P}}
\newcommand{\QQ}{{\cal Q}}
\newcommand{\RR}{{\cal R}}
\newcommand{\SS}{{\cal S}}
\newcommand{\TT}{{\cal T}}
\newcommand{\TTT}{{\cal P}}
\newcommand{\UU}{{\cal U}}
\newcommand{\VV}{{\cal V}}
\newcommand{\WW}{{\cal W}}
\newcommand{\XX}{{\cal X}}
\newcommand{\YY}{{\cal Y}}
\newcommand{\ZZ}{{\cal Z}}

\newcommand{\A}{{\Bbb A}}
\newcommand{\C}{{\Bbb C}}
\newcommand{\D}{{\Bbb D}}
\newcommand{\Hyp}{{\Bbb H}}
\newcommand{\J}{{\Bbb J}}
\newcommand{\Ll}{{\Bbb L}}
\renewcommand{\L}{{\Bbb L}}
\newcommand{\M}{{\Bbb M}}
\newcommand{\N}{{\Bbb N}}
\newcommand{\Q}{{\Bbb Q}}
\newcommand{\R}{{\Bbb R}}
\newcommand{\T}{{\Bbb T}}
\newcommand{\V}{{\Bbb V}}
\newcommand{\U}{{\Bbb U}}
\newcommand{\W}{{\Bbb W}}
\newcommand{\X}{{\Bbb X}}
\newcommand{\Z}{{\Bbb Z}}

\newcommand{\tT}{{\mathrm{T}}}
\newcommand{\tD}{{D}}
\newcommand{\hyp}{{\mathrm{hyp}}}

\newcommand{\f}{{\bf f}}
\newcommand{\g}{{\bf g}}
\newcommand{\h}{{\bf h}}
\renewcommand{\i}{{\bar i}}
\renewcommand{\j}{{\bar j}}

\def\Bf{{\bold{f}}}
\def\Bg{{\bold{g}}}
\def\BG{{\bold{G}}}
\def\Bh{{\bold{h}}}
\def\BT{{\bold{T}}}
\def\Bj{{\bold{j}}}
\def\Bphi{{\bold{\Phi}}}
\def\Bpsi{{\bold{\Psi}}}
\def\B0{{\bold{0}}}

\newcommand{\Com}{\cal Com}
\newcommand{\Top}{\cal Top}
\newcommand{\QC}{\cal QC}
\newcommand{\Def}{\cal Def}
\newcommand{\Teich}{\cal Teich}
\newcommand{\QL}{{\cal{QG}}}
\newcommand{\PPL}{{\cal P}{\cal L}}

\newcommand{\hf}{{\hat f}}
\newcommand{\hg}{f}
\newcommand{\hz}{{\hat z}}
\newcommand{\hM}{{\hat M}} 
\newcommand{\hV}{{U}}
\newcommand{\hVV}{{\cal U}}

\renewcommand{\lq}{``}
\renewcommand{\rq}{''}


\catcode`\@=12

\def\Empty{}
\newcommand\oplabel[1]{
  \def\OpArg{#1} \ifx \OpArg\Empty {} \else
  	\label{#1}
  \fi}
		
%

\long\def\realfig#1#2#3#4{
\begin{figure}[htbp]
\centerline{\psfig{figure=#2,width=#4}}
\caption[#1]{#3}
\oplabel{#1}
\end{figure}}

%

\newcommand{\comm}[1]{}
\newcommand{\comment}[1]{}

\bigskip\bigskip

\title[Parabolic renormalizations]{Note on the geometry of generalized parabolic towers}
\author {Mikhail Lyubich }
\date{January 5, 2002}

\begin{abstract} In this technical note we show that  the geometry of 
generalized parabolic towers cannot be essentially bounded. 
It fills a gap in \cite{attrac}. 
\end{abstract}

\maketitle


\section{Introduction}

In 1985 Milnor posed a problem of existence of ``wild attractors'' for real quadratic polynomials
\cite{M}. In this situation  a quadratic polynomial $f$ would have a topologically transitive  invariant
interval $I$ such that almost all orbits on $I$ would  converge to a Cantor set $A\subset I$. 
It was proven in \cite{attrac}
that this is impossible (Theorem II). This result was deduced from a fundamental property
of {\it exponential decay of geometry} (Theorem I)  which has found many further applications since then.    

The above results were stated in \cite{attrac} for a more general class of maps, 
 $S$-unimodal maps with  non-degenerate critical point. The $S$-unimodal case  was reduced to the  quadratic one
by taking (generalized) quadratic-like limits of (generalized) renormalizations and 
applying to them a ``quasiconformal trick''.   
In this reduction one special case   
when the limits of generalized renormalizations become parabolic
was missed from the consideration 
(this gap was pointed out by Oleg Kozlovski). In this note we will show how to fill this gap. 

To this end we replace the generalized renormalization with 
   the "generalized parabolic renormalization"  
relating different parabolic  maps obtained in the limit. 
In this way we obtain a ``parabolic tower'' of generalized  quadratic-like  maps.
Applying  the quasi-conformal trick to such a tower, 
we conclude that it has  exponentially decaying geometry.
(The idea of the quasiconformal trick  is  to show that any two towers with the same combinatorics are 
quasiconformally equivalent and then to provide one example of a tower with a
given combinatorics and decaying geometry.)

Note that the part of the argument which contains the
gap 
 is concerned with the smooth case only.
In the case of quadratic polynomials (or, more generally, analytic maps of Epstein class),
the proof of \cite{attrac} does not need any adjustments.

In \S \ref{outline} we will outline, for reader's convenience,  the original argument of \cite{attrac}. 
In the following two sections, \S\ref{parabolic towers} and \S \ref{complex towers}, we will review the theory
of generalized parabolic towers, real and complex. In \S \ref{trick} we will prove, by means of the quasiconformal  trick, that such towers
have decaying geometry. In the Appendix we will  give another proof of this result,  
by adjusting complex estimates from \cite[\S 6]{puzzle}. 

\msk
{\it Remark.}  Yet another way to prove decay of geometry of $S$-unimodal maps
is by means of the asymptotically
conformal extension of  smooth maps as proposed in \cite[\S 12.2]{puzzle}.

\medskip {\bf Acknowledgement.} I thank the participants of the seminar
in IMPA (August 2000), particularly A.~Avila, O.~Kozlovski, W.~de~Melo, and M.~Shishikura,
for their attention and valuable comments.    
 
\section{Outline of the original argument}\label{outline}
Let $f: [-1,1]\ra [-1,1]$ be a  $C^3$ unimodal map  with negative Schwarzian derivative and  non-degenerate critical point at 0
such that $f(1)=1$ and $f'(1)>1$. 
Assume for simplicity that $f$  is even.
Such a map has a unique fixed point $\alpha\in (-1,1)$.
 Assume this point is also repelling and 
$f$ is non-renormalizable with recurrent critical point. 

Let $I^0=[\alpha, -\alpha]$. Define inductively the {\it principal nest} of $0$-symmetric intervals
$$  I^0\supset I^1\supset I^2\supset\dots  $$
as follows. Let $r_n$ be the first landing time of the critical point at $I^{n-1}$.
Then $I^n$ is the component of $f^{-r_n} (I^{n-1})$ containing $0$. Let $g_n= f^{r_n} | I^n$. 

There are  two  different types of returns of the critical point to $I^{n-1}$, central and non-central.
The return to $I^{n-1}$ (and the corresponding level $n-1$) is called {\it central } if  $g_n(0)\in I^n$. 
Let $\{ n_k - 1\} $ be the sequence of non-central levels in the principal nest. 

The ratios $\la_n= |I^n| / |I^{n-1}|$ are called the {\it scaling factors} of $f$. 

\begin{thm}[Theorem II of \cite{attrac}]\label{II}
  There exist constants $\rho\in (0,1)$ and  $C>0$
such that  $\la_{n_k+1}\leq C\rho^k $.
\end{thm}
 
We say that maps satisfying the conclusion of Theorem II have {\it (exponentially) decaying geometry}. 

The proof of this theorem given in \cite{attrac} contained a gap which will be explained and filled below. 

Theorem II was first proven under the assumption that one of the scaling factors
is sufficiently small:

\begin{lem}\label{conditional II}
   There exists an absolute constant $\de>0$ 
such that if $\la_N < \de$ for some $N$  then  $f$ has exponentially decaying geometry.
\end{lem}  

The further analysis depended on the combinatorics of $f$ defined in terms of generalized renormalizations
$g_n : \cup I^n_i \ra I^{n-1}$ of $f$. 
Consider the sequence of first return maps $f_n : \cup I^n_i \ra I^{n-1}$ to the intervals of the principal nest.
For each $n$,  $I^n_i$ are disjoint closed intervals,  
$I^n_0\equiv I^n$, the map $f_n|I^n_0$ coincides with the previously introduced $g_n$,
and $f_n$ diffeomorphically maps each non-central interval $I^n_i$, $i\not=0$,  onto $I^{n-1}$. 
The {\it generalized renormalization} $g_n$ is the restriction of $f_n$ to the union of intervals $I^n_i$ intersecting
the critical orbit. 

Let us say that $f$ has {\it essentially bounded geometry} if all the non-central intervals $I^n_i$, $i\not=0$,
 and all the gaps
in between (i.e., the components of $I^{n-1}\sm \cup I^n_i$) are commensurable (uniformly on all levels). 
Note that in the case of essentially bounded geometry,  there are only finitely many intervals $I^n_i$ on each level
(and their number is bounded).

\begin{lem}[\cite{attrac}, \S 3]\label{unbounded geom}
  Maps with essentially unbounded geometry have exponentially decaying geometry.  
\end{lem}

To treat the bounded case we passed to limits of the renormalizations.
These limits have better qualities than the original maps. 

A map $g: \cup_{i=0}^l V_i\ra \De$ is called a (real)
{\it generalized quadratic-like} map (see \cite{puzzle})  if $V_i$ and $\De$ are ($\R$-symmetric) topological disks in $\C$,
$V_i\Subset \De$, $\bar V_i$ are pairwise disjoint,
$g: V_0\ra \De$ is a double branched covering, while the maps $g: V_i\ra \De$ 
are conformal isomorphisms (preserving the real line). 
We put the critical point of $g|V_0$ at the origin.  In what follows, all generalized quadratic-like maps are assumed to be real.

\begin{lem}[\cite{attrac}, \S\S 3 - 4]\label{renorm limits}
  If $f$ has essentially bounded geometry then the renormalizations $g_n$ have an analytic limit which admits 
an extension to a generalized quadratic-like map $g: \cup V^n_i \ra \De$.
\end{lem}

If the limit $g$  of renormalizations has decaying geometry then the original map $f$ has arbitrary small scaling factors
and hence has exponentially decaying geometry by Lemma \ref{conditional II}. 
This reduces Theorem \ref{II} to the quadratic-like case.
We treated this case by means of the ``quasiconformal trick'' based on the following lemma:

\begin{lem}[Kahn, see Lemma 5.3 of \cite{attrac}]\label{Kahn}
  Any two non-renormalizable generalized quadratic-like maps  $g$ and $\tl g$ with the same combinatorics 
are qc  conjugate. \footnote{``qc'' stands for ``quasiconformal''.}  
\end{lem}  

Since the property of  decaying geometry is qc invariant, it is enough to give one example of a quadratic-like map with 
giving combinatorics and decaying geometry to conclude that this property is valid for all maps with that combinatorics. 
But it is easy to show (using a version of the kneading theory, see \cite{attrac}, \S 1) that for any given combinatorics,
there exists a generalized quadratic-like map $g: \cup V^n_i\ra \De$ with this combinatorics
whose  central branch $g|V^n_0$ is purely quadratic while all non-central branches $g|V^n_i$ are linear. 
Moreover, the central domain $V^n_0$ can
be selected arbitrary small in $\De$, so that $g$ has exponentially decaying geometry by Lemma \ref{conditional II}. 

\msk
What is overlooked in the above argument is that the limit quadratic-like map $g$ in Lemma \ref{renorm limits}
 does not have to be non-renormalizable.
However, as the following simple lemma shows (compare [\cite{puzzle}, Theorem V]), the only alternative is parabolic:

\begin{lem}\label{parabolic cascades}
Assume $f$ has essentially bounded geometry. 
Let $n-1$ be a central level of the principal nest and $m-1$ be the following non-central level, $m>n+N$.
If $N$ is sufficiently big   then $0\not\in g_n(I^n)$ (``low return''), and
 $g_n | I^n $ is close to a unimodal map with parabolic fixed point. 
\end{lem}   
 
Passing to limits of generalized renormalizations on consecutive non-central levels,
we  obtain a ``parabolic tower'' consisting of maps related by either the above generalized renormalization or
by its parabolic version. 
Applying the quasiconformal trick to such towers, we complete the proof of Theorem \ref{II}. 

Below we will elaborate on  this argument. 

\comm{
We say that combinatorics of  $f$  is {\it bounded by } $N$   if:

\ssk\nin $\bullet$
  The number of intervals $I^n_i$ in the $\Dom(g_n)$ is bounded by $N$;

\ssk\nin $\bullet$
  The return times of intervals $I^{n+1}_j$ back to $I^n$ under the iterates of $g_n$ are bounded by $N$.
}

\medskip
In the case when $f$ is a quadratic polynomial  (or more generally, an analytic map  of Epstein class),
Lemma \ref{renorm limits} is not needed since 
$f$ itself has a generalized  quadratic-like  renormalization $g$ to which we can apply the quasiconformal trick. 
So, in this case, the proof of Theorem II given in \cite{attrac} was complete.

\section{Real parabolic towers}\label{parabolic towers}

\subsection{Definition}

Let $\RR$ be the class of $C^3$-smooth one-dimensional maps (considered up to rescaling)
$$
    g: \bigcup_{i=0}^{l} I_i \ra T, 
$$
with negative  Schwarzian  derivative such that

\ssk\nin $\bullet$
the $I_i$ are disjoint closed intervals contained in $\inter T^{n}$; $I_0\ni 0$;

\nin $\bullet$
the restrictions  $g : I_i\ra T$ are diffeomorphisms for $i\not=0$;

\nin $\bullet$ 
the restriction $g : (I_0, \di I_0) \ra (T, \di T)$ is an even unimodal map with a
{\it non-degenerate}  critical point at 0;

We let $\RR(\eps)$ be the subclass of  maps which admit an extension $\hat I_i\ra \hat T$
of class $\RR$ with  $\hat T = (1+\eps) T$ ($\eps$ is called the ``extension parameter''). 

\msk

Consider a nest of closed intervals $T^n$ containing 0, 
and  a sequence of maps  
$$
    g_n: \bigcup_{i=0}^{l_n}  I^n_i  \ra T^{n}, 
$$
of class $\RR(\eps) $ for some $\eps>0$ (independent of the level). 
We assume that  the next map, $g_{n+1}$, is related to the previous one, $g_{n}$,
by a generalized renormalization of three possible types,
standard, cascade, or parabolic. We will now describe these renormalizations. 

\ssk
The {\it standard} renormalization is applied in the non-central case, i.e., when $g_n(0)\not\in I^n_0$. 
Then  $T^{n+1} = I^{n}_0$,
 and $g_{n+1} $ is  the first return map to $T^{n+1}$  under iterates of $g_n$
restricted to some 
intervals $I^{n+1}_i$.

\ssk
The  {\it cascade} renormalization is applied in the central escaping  case 
i.e., when $g_n(0)\in I^n_0$ but 
there is some $t>0$ and $j\not=0$ such that   $ (g_n\,|\, I^n_0)^{\circ t} (0) \in I^n_j.$
In this case, $T^{n+1}$  is the pull-back of $I^n_j$ under $g_n^{\circ t}$ containing 0 
(i.e., the component of $g_n^{-t} I^n_j$ containing 0),
and $g_{n+1}|I^{n+1}_i$ is composed of some 
branches of the first return map to $T^{n+1}$.

In this case  we also consider  a bigger family $\LL^{n}$ of intervals
$L^{n}_s$ obtained from the  non-central intervals $I^{n}_j$, $j\not=0$, as the pull-backs 
under  iterates of $g_n| I^{n}_0$. 
We define the {\it transit map }  as follows: $\Psi_{n}| L^{n}_s = g_n^{\circ t} | L^{n}_s$
where $g_n^{\circ t}$ maps $L^{n}_s$ into  $I^{n}_j$.
One of these intervals coincides with $T^{n+1}$; it will be denoted $L^{n}_0$. 
The times $t=t(n,s)$ will be called the {\it transit times} of the intervals $L^n_s$.

Let $G_{n} = g_{n}\circ \Psi_{n}$. This is a ``Bernoulli map'' which diffeomorphically 
maps each non-central interval $L^{n}_j$ onto $T^{n}$ and unimodally maps
$(L^{n}_0, \di L^{n}_0)$ into $(T^{n}, \di T^{n})$.  

\ssk
The {\it parabolic} renormalization is applied  when $g_{n}: I^{n}_0\ra T^{n}$
 is a parabolic map supplied with a Lavours  map $\tau_n$ through  the parabolic point
(see \cite{D1,Sh}). In this case we consider pull-backs $L^{n}_s$ of the
non-central intervals $I^{n}_j$ under joint iterates of $g_n |I^n_0$ and
$\tau_n$. Thus, for any $L=L^n_s$, there exist $k=k$ and $l$ such that the transit map
$$ 
\Psi_n|L \equiv g_n^{\circ l} \circ \tau_n  \circ     g_n^{\circ k} |\; L $$
is either a diffeomorphism onto some $I^n_j$ (for $s\not=0$) 
or  is a unimodal map (for $s=0$).
Let $T^{n+1} = L^{n}_0$ and let $G_{n}= g_{n}\circ \Psi_{n}$  be the associated  Bernoulli map. 
Then $g_{n+1}$  is composed of some branches of the first return map to $T^{n+1}$ 
under iterates of $G_{n}$. 
   
Let us call such a sequence of maps $(g_n, \Psi_n)$ a (one-sided {\it parabolic) tower} $\GG$
of class $\TT(\eps)$.  
(We will usually skip the transit maps  $\Psi_n$ from the notation.) Let $\TT=\cup \TT(\eps)$.
In the case when there are no parabolic maps among the $g_n$, we we also refer to the tower as 
a {\it principal nest} (of generalized renormalizations of the top map).

We will consider towers up to rescaling of the base interval $T^0$.

Let $\II^n$ stand for the family of intervals $I^n_i$ and let $\II^n_*= \II^n\sm \{I^n_0\}$. 
In the standard case, 
we let $G_{n} \equiv g_{n}$ and $L^{n}_j\equiv I^{n}_j$. 

\subsection{Combinatorics of towers}
The combinatorics $\kappa(g_n)$  is determined by the 
order of the intervals $L^{n}_s$ on the real line and by the
itineraries of the intervals $I^{n+1}_i$ through the  $L^{n}_s$
under the iterates of $G_{n}$ until the first return back to $T^{n+1}$.
The combinatorics $\kappa= \kappa(\GG)$ of a tower $\GG=\{g_n\}$  is the sequence 
$\{ \kappa(g_n) \}$. 

Given an interval $L\in \LL^{n} $, 
let $k_+$ be defined by the property  that $g_{n}^{\circ k_+} L \in \II_*$ 
(if such $k_+$ does not exist, we let $k_+=\infty$),
 and let $k_-$ be the biggest $k$ such that
$L = g_{n}^{\circ k} K$ for some interval $K\in \LL^{n}$. 
We define the {\it depth} of $L$ (and of any point $x\in L$) as $\min(k_-, k_+)$.     

Two combinatorics on a given level are called {\it essentially equivalent} if they become the same
after removing all intervals of sufficiently big depth. 
Two tower combinatorics are essentially equivalent if they are such on every level.

  Let us consider a combinatorics $\kappa$ of a parabolic tower and a sequence 
of essentially equivalent combinatorics $\kappa_n$ of principal nests. We say that
$\kappa_n\to \kappa$ if the transit times of the $\kappa_n$ on the parabolic levels 
of $\kappa$ grow to $\infty$. 
(See \S 3 of \cite{Hi} for a more detailed discussion of the space of combinatorics.)

\ssk
The combinatorics of a tower is called {\it essentially bounded} if:

\nin$\bullet$ The numbers $l_n$ of the intervals $I^n_i$ are bounded;

\nin$\bullet$ 
 If an interval  $J\subset I^n_j$ belongs to the itinerary of some $I\in \II^{n+1}$
     then  $g_n J$  lands at an interval $L\in \LL^n$ of bounded depth;

\nin$\bullet$
  The return times of the intervals $I\in \II^{n+1}$ back to $T^{n+1}$ under the iterates of  $G_n$ are bounded. 

If all the above numbers are bounded by $p$ we say that the essential combinatorics
of the tower is bounded by $p$.

\subsection{Essentially bounded geometry}

We say that a parabolic tower has an {\it essentially bounded geometry} if: 

\msk
\noindent $\bullet$ the intervals $I^n_j$ and the gaps in between them 
(i.e., connected components of $T^n\sm \cup I^n_j$)
are all commensurable with $T^n$;

\noindent $\bullet$
 all the maps $(g_n|I^n_j)^{-1}$ ($j\not=0$), $\tau_n^{-1}$ admit a diffeomorphic
 extension with positive Schwarzian derivative to 
a definitely bigger neighborhoods of their domains of definition;

\noindent $\bullet$ $g_n | I^n_0 $ is a composition of the quadratic map and a diffeomorphism $h_n$
  whose inverse admits an extension as above.  b
\msk

By an obvious quantification of  this notion (which also incorporates the extension parameter $\eps$ of the tower),
 we can make sense of  the notion of ``essentially $C$-bounded geometry''
Moreover, this notion makes sense on any particular  level of the tower.  
Let $\TT(p,C)$ stand for the space  of parabolic towers
with  essentially $C$-bounded geometry on the top level
and essentially $p$-bounded combinatorics. 
The following statement is a slightly modified Lemma 8.8 of \cite{puzzle}.

\begin{lem}\label{ess bounds}
Given $C$ and $p$, towers $\GG\in \TT(p,C)$ 
have essentially bounded geometry on every level. 
\end{lem}

 We are ready to formulate the main result of this note:

\begin{thm}\label{filling}
  There are no parabolic towers with essentially bounded geometry.
\end{thm}

This theorem together with Lemma \ref{unbounded geom} imply Theorem \ref{II}.

\subsection{Scaling factors}

The ratios $\la_n= |T^n| / |T^{n-1}|$ are called the {\it scaling factors} of 
the tower. 

\begin{lem}\label{conditional result}
   There exists an absolute constant  $\de>0$ 
such that if $\la_1< \de$ then  $\la_n\leq C\rho^n $ for some $\rho\in (0,1)$, $C>0$.
\end{lem}

This is a straightforward generalization of Lemma \ref{conditional II}. 
Note that in the cascade case, the proof  of Lemma \ref{conditional II}  uses only the transit map
through the cascade, so that it applies to the parabolic case as well.

\section{Complex parabolic towers}\label{complex towers}

\subsection{Analytic limits of Epstein class}
A sequence of towers $\GG_m=\{g_{m,n}, \Psi_{m,n}\}$ converges to a tower $\GG=\{g_n, \Psi_n \}$ if on every level $n$
the domains of the maps $g_{m,n}$ converge to the domain $D_n$  of $g_n$ in the Hausdorff topology,
$g_{m,n}\to g_n$  in the $C^1$-topology on the components of $D_n$, 
and $\Psi_{m,n}\to \Psi_n$ in the $C^1$-topology on the components of the domain of $\Psi_n$.  


\ssk
A parabolic tower is called a tower of  {\it Epstein class} if

\ssk\nin $\bullet$  all the maps $g_n$ and  $\Psi_n$ are analytic;

\ssk\nin $\bullet$ the inverse maps $( g_n | I^n_j)^{-1}$, $j\not=0$,
  admit a holomorphic extension to the slit plane $\C\sm (1+\eps)T^n$,
  and the similar statement holds for the transit maps $\Psi_n$ ;

\ssk\nin $\bullet$ the map $g_n | I^n_0$ is a composition of the quadratic map and a diffeomorphism 
     which admits a holomorphic extension as in the previous item.
 
By the $m$-shifted tower $\GG^m$ we mean the tower $\{ g_n\}_{n\geq m}$.
Now, the Shuffling Lemma (see \cite[Ch. VI, Theorem 2.3]{MvS}) yields:

\begin{lem}\label{Epstein limits}  
Let $\GG$ be a parabolic tower with bounded geometry. 
Then any sequence of shifted towers $\GG^m$ contains a subsequence converging to a tower of Epstein class.
\end{lem}

\subsection{Generalized quadratic-like maps}\label{gpr}

We will now complexify the discussion of \S \ref{parabolic towers}.

Let $\CC$ be the space of (real symmetric) generalized quadratic-like maps (up to rescaling of $\De$).
Let $\CC(\eps)$ be the subspace of  maps which admit an extension $\hat V_i\ra \hat \De$
of class $\CC$ with $\mod(\hat \De \sm \De)>\eps$.

This notion immediately allows us to complexify the notion of parabolic tower.
Such a tower  consists of generalized quadratic-like maps $g_n: \cup V^n_i\ra \De^n$ related by either standard, or
cascade, or parabolic renormalization, and associated transit maps $\Psi_n: U^n_s \ra L^n_j$,
where $U^n_s$ are topological disks whose real slices coincide with the intervals $L^n_s$. 
We will also consider the complexification of the associated Bernoulli maps $G_n: \cup U^n_s\ra \De^n$.

\comm{
  Consider a generalized quadratic-like map $ g: \cup V_i \ra \De$  such that
the restriction $g:  V_0\ra \De$ is parabolic. 
Let $A= \De\sm \bar V_0$. We denote
by $\VV$ the family of domains $V_i$ and by $\VV_*$ the family of off-critical domains
$V_i$, $i\not=0$. Let us fix a transit map $\tau$ through the
parabolic point $\alpha$. 

Let us also consider a bigger family  $\hVV\supset \VV$ of topological disks
$\hV_j\Subset \De$  and
 a map $\hg: \cup \hV_j\ra \De$ satisfying the following properties (see Figure 1):

\ssk
 
\nin $\bullet$  the map $\hg |\, \hV_j$ is a restriction of either $g$ or $\tau$;

 \nin $\bullet$  $\hg$ univalently maps $\hV_j\in \hVV_*\sm \VV$   onto some $\hV_{k(j)}$
  (where $\hVV_*$ stands for the family of off-critical domains of $\hVV$);

 \nin $\bullet$  $\hg | \hV_0$ is a double branched covering over some $\hV_{k(0)}$;
  
\nin $\bullet$ 
   for any $\hV_j\in \hVV$ there is a time $t(j)\leq N$ such that 
           $\hg^{t(j)} \hV_j \in \VV$.

The last property allows us to  define the transit map $\Psi: \cup U_j\ra \cup V_i$,
 $\Psi | U_j = \hg^{t(j)}$.

 Let us then consider a Bernoulli map $G: \cup \hV_j \ra \De$ defined as follows:
$G = g\circ \Psi$. It univalently maps each off-critical domain 
$\hV_j\in \hVV$ onto $\De$, and it is a double branched covering of $\hV_0$ over
$\De$. 

Let $L: \cup \Om_k \ra \hV_0$ be the first landing map at $\hV_0$, i.e., 
$L(z) = G^{l(z)} (z)$, where $l(z)$ is the first time (if it exists)
when the $G$-orbit of $z$ lands at the critical domain $\hV_0$. It is defined on the
union of disjoint topological disks $\Om_k\subset \cup \hV_j$.  


We can now consider the {\it first return map} $R: \cup V_s'\ra \De'\equiv \hV_0$
as the composition $L \circ G$. Here the domains $V_s'$ are the pull-backs of the
$\Om_k$ under $G| \hV_0$. Finally, let us restrict $R$ to the union of domains $V_s'$
that intersect  the critical orbit $\{R^n 0\}$. We obtain a generalized parabolic 
renormalization $g':  \cup V_s'\ra \De'$ of $g$. Denote by $\VV'$ the family of domains
$V_s'$ on which  $g'$ is defined.  
}

\subsection{A priori bounds and rigidity}

By repeating the argument of \cite[\S 4]{attrac} 
(using on parabolic levels the parabolic transit maps instead of almost parabolic transit maps), we obtain:

\begin{lem}[Complex bounds]\label{complex bounds}
  If $\{ g_n\}$ is a parabolic tower of Epstein class then one of the maps
$g_n$ admits a  generalized quadratic-like extension.
\end{lem} 

{\it Remark.} One can also use complex bounds  of \cite{LY} to prove the above statement.

\begin{prop}[\cite{Hi}, Prop. 6.1]\label{rigidity}
  A complex parabolic tower is uniquely determined by its combinatorics 
and the top map $g_1$.
\end{prop}

\subsection{Space of towers}

Similarly to the real case, 
the space of generalized quadratic-like maps $g_m: \cup_{i=0}^l V_i\ra \De$ 
and the space of complex towers are endowed with  natural topologies. 

Given a generalized quadratic-like map $g: \cup V_i\ra \De$, 
let $A_i\Subset \De\sm \overline{\cup V_i}$ be an annulus of maximum modulus 
surrounding $V_i$ but not surrounding the other domains $V_j$, $j\not=i$.
Let $\mod(g)\ \min \mod (A_i)$. 

\begin{lem}\label{compactness}
  The space of complex parabolic towers $\GG\in \TT(p,C)$
 with a definite modulus on the top level is compact.
\end{lem} 

\begin{pf}
A simple estimate shows that $\mod (g_{n+1})\geq {1/2} \mod(g_n)$
so that the moduli of the tower maps on every  level are definite.
Moreover, the extensions for the inverse maps $(g_n | V^n_i)^{-1}$ provide us
with  extensions of the transit maps with a definite modulus. 
Now the assertion follows from the  Carath\'eodory  compactness of 
the family of  pointed domains with bounded from below inner radius (see \cite[Thm 5.2]{McM})
and  normality argument. 
\end{pf}

Let us consider a holomorphic family of generalized quadratic-like maps 
$g_\la: \cup V_{i,\la}\ra \De_\la$ over a topological disk $\La\subset \C$
such that the boundary  $\di \De_\la$ moves holomorphically with $\la$.  
Such a family is called {\it proper} 
if $g_\la(0)\to \di \De_\la$ as $\la\to \di \La$. A proper family is called
{\it unfolded} if for a Jordan curve $\gamma\subset \La$ close to $\di \La$,
the curve $\la\mapsto g_\la(0)$, $\la\in \gamma$, has winding number 1
around 0 (see \cite{DH,parapuzzle}). 

\begin{lem}\label{approx}
  Let $\GG$ be a  complex parabolic tower with combinatorics $\kappa$.
Let  $\{\kappa_l\}$ be a sequence of essentially equivalent nest combinatorics 
converging to $\kappa$. 
If these combinatorics have a sufficiently big transit time on the top parabolic level,
then there exists a sequence $\GG_l$ of complex principal nests with combinatorics 
$\kappa_l$ converging to $\GG$.
\end{lem}

\begin{pf}
 Let $g_s$ be  the  first parabolic map in the tower. 
 Include the map $g_1\equiv g_{1,0}$  into a one-parameter family of generalized quadratic-like maps $g_{1,\la}$
over some domain $\La$ with efficiently changing  multiplier of the parabolic point of $g_s $
(i.e., the parameter derivative of this multiplier does not vanish).  
Then there is a sequence of domains $\Om_m\subset \La$ converging to 0 
such that the combinatorics of the tower does not change on levels $1,2,\dots, s$,
 and the associated  Bernoulli maps $F_{m,\la}$ on level $s$  form proper unfolded families  over the $\Om_m$ (compare \cite{D2}).
These domains are  labelled by the transit time $m$ of the critical point through 
the almost parabolic region until it lands at the appropriate non-central  domain of $g_{s,\la}$.  
If the transit time of the combinatorics $\kappa_l$ on  level $s$  is sufficiently big,
then we can select $m$ to be equal to this transit time.

Restricting the family $F_{m,\la}$ to the real slice of $\Om_m$, we obtain a full family of 
maps of class $\RR$. 
Since any admissible nest combinatorics is realizable in such a family  (see \cite[\S 1]{attrac}),
the  desired combinatorics $\kappa_l$ are all realizable by some principle nests $\GG_l$. . 

 By Lemma \ref{compactness}, the sequence $\{\GG_l\}$  is pre-compact.  
Moreover, any limit parabolic tower has combinatorics $\kappa$ 
and by Proposition \ref{rigidity} coincides with $\GG$. 
Hence   $\GG_l\to \GG$.  
\end{pf}

\section{Quasi-conformal trick}\label{trick}

We will now adjust the argument of \cite[\S 5]{attrac} to the case of parabolic towers. 

A homeomorphism  $h: \De\ra \tl\De$ is called a {\it pseudo-conjugacy} between two generalized
quadratic-like maps $g: \cup V_i \ra \De$ and $\tl g: \cup \tl V_i \ra \tl\De$
if $h$ is equivariant on the boundary, i.e., $h(gz) = \tl g(hz)$ for $z\in \cup \di V_i$.

\msk{\it Remark.} 
 Under the circumstances of Lemma \ref{Kahn}, any qc pseudo-conjugacy between $g_0$ and $\tl g_0$
promotes to a qc conjugacy with the same dilatation.

\begin{lem}\label{qc conjugacy 1}
  Any two parabolic towers with the same combinatorics are qc conjugate.
\end{lem}\label{parab conjugacy}

\begin{pf} Consider two parabolic towers $\GG$ and $\tl \GG$ with the same combinatorics. 
By Lemma \ref{approx}, they can be approximated by combinatorially equivalent principal nests
$\GG_n$ and $\tl \GG_n$. By the above Remark, these nests are qc equivalent
with bounded dilatation. Passing to limits, we conclude that the towers $\GG$ and $\tl \GG$
are qc equivalent as well. 
\end{pf}

{\it Proof of Theorem \ref{filling}}.
If there exists a real parabolic tower with essentially bounded geometry, 
then by Lemmas \ref{Epstein limits} and \ref{complex bounds}, 
there exists a complex parabolic tower with this property as well. 
We will now show by means of the quasi-conformal trick that such complex towers do not exist. 

 It is easy to create a complex parabolic tower $G$ with a 
given combinatorics and an arbitrary small initial scaling factor $\la_1$.    
Namely, consider a generalized quadratic-like map $g: \cup V_i \ra \De$ with 
the given combinatorics whose central branch $g| \De$ is   purely quadratic,
while the non-central  branches $g| \V_i$, $i\not=0$, are linear. 
The central domain $V_0$ of this map can be selected as an arbitrary small round disk.    
Changing the transit map, we can produce a full family of towers 
on the consecutive levels (compare the proof of Lemma \ref{approx}), 
so that there is a tower with the required combinatorics among them. 

By Lemma \ref{conditional result}, the scaling factors of $\GG$ decay exponentially. 

Let $\tl \GG$ be an arbitrary parabolic tower with the same combinatorics. 
By Lemma \ref{qc conjugacy 1}, $\tl \GG$ is qc equivalent to $\GG$.
Since the exponential decay of the scaling factors is a qc invariant property,
it holds for $\tl \GG$ as well. 
\QED

\section{Appendix} 

\subsection{Isles and asymmetric moduli}
Here we reproduce for convenience some definitions from \S 5
of  \cite{puzzle}.

 Let $\{V_i\}_{i\in \II}\subset \VV$  be a finite family of disjoint 
 domains consisting of at least two elements (that is $|\II|\geq 2$)
and containing a critical domain $V_0$.
Let us call such a family {\it admissible}.
We will freely identify the label set $\II$ with the family itself.

Given a domain $D\subset \De$ whose boundary does not intersect $V_i$, 
let $\II|D$ denote the family of domains of $\II$ contained in $D$.
Let $D$  contain at least two domains of family $\II$.
For $V_i\subset D$,   let 
$$R_i\equiv R_i(\II|D)\subset D\backslash\bigcup _{j\in \II|D}V_j$$ be
an annulus of maximal modulus enclosing $V_i$  but not enclosing other
domains of the family $\II$. We will briefly call it the 
{\sl maximal annulus} enclosing
$V_i$ in $D$ (rel the family $\II$).

Let us define
{\sl the asymmetric modulus of the family $\II$ in $D$ } as
$$\sigma(\II|D)=\sum_{i\in \II} {1\over2^{1-\delta_{i0}}}
  \mod R_i(\II|D),$$      
where $\delta_{ji}$ is the Kronecker symbol. So  the critical
modulus is supplied with weight 1, 
while the off-critical moduli are supplied with weights 1/2
(if $D$ is  off-critical then all the weights are
actually 1/2). 


For  $D=\De$,  let $\sigma(\II)\equiv\sigma(\II| \De )$.
The asymmetric modulus of $\VV$ is defined as follows: 
$$\bar\sigma(\VV) = \min_\II\sigma(\II),$$     
where $\II$ runs over all admissible sub-families of $\VV$.

Let $\{V'_i\}_{i\in \II'}$ be an admissible sub-family of $\VV'$.
Let us organize the domains of this family in {\sl isles}
in the following way. 
A domain $D'\subset \Delta'$ is called an {\it island} (for family
$I'$) if

\smallskip\noindent $\bullet$ $D'$ contains at least two domains of
  family $\II'$;

\smallskip\noindent $\bullet$ There is a $t\geq 1$ such that $G^k D'\subset
  V_{i(k)}, \; k=1,\ldots t-1,$ with $i(k)\not=0$, while $G^t D=\Delta$
  (where $G$ is the Bernoulli map from \S \ref{gpr}).   

Given an island $D'$, let 
$ \phi_{D'} = G^t : D'\rightarrow\Delta$. This map
is either a double covering or a biholomorphic isomorphism
 depending on whether $D'$ is critical or not.  In the former case,
 $D'\supset V'_0$. 

We call a domain $V_j'\subset D'$  $\phi_{D'}$-{\sl precritical} 
if $\phi_{D'}(V_j')=V_0$. 
There are at most two precritical domains in any $D'$.
If there are actually two of them, then they are off-critical and
symmetric with respect to the critical point 0. In this case
$D'$ must  also contain the critical domain $V_0'$.


Let ${\cal D}'={\cal D}(\II')$ be the family of isles associated with $\II'$.
Let us  consider the asymmetric moduli $\sigma(\II'|D')$ as a function on
this family. This function is clearly monotone:
\begin{equation}\label{5-3}
\sigma(\II'|D')\geq\sigma(\II'|D_1')\quad \text{if} \quad D'\supset D_1',
\end{equation}

Let us call an island $D'$ {\it innermost} if it does not contain any
other isles of the family $\DD(\II')$. As this family is finite,
innermost isles exist.

\subsection{Growth  of moduli}
Let $\mu=\mod(A)\equiv \mod(\De\sm V_0)$, $\mu'=\mod(A')\equiv \mod(\De'\sm V_0')$,
$\nu= \mod(\De \sm U_0)$.

\begin{lem}\label{mu}
$  \mu' \geq {1\over 2}\nu \geq \mu. $
\end{lem}

\begin{pf}
Let $B$ be the parabolic basin of $g_0\equiv g|V_0$.
Then 
 by the Gr\"otzsch inequality, 
$$\nu \geq \mod(\De\sm B) \geq \sum_{n=0}^\infty \mod (g_0^{-n} A) = 2\mu, $$ 
which is the second required inequality.

To prove the first one, consider a double covering  $G: \De'\sm V_0' \ra  \De\sm \Om_k$
(where $\Om_k$ is an appropriate domain of the first landing map $L$ from  \S \ref{gpr}). 
Since $\mod(\De\sm \Om_k)\geq \nu$, the first inequality follows as well.
\end{pf}
  
The following statement is Lemma 5.4 from \cite{puzzle}.
 
\begin{lem}\label{main} Let $\II'$ be an admissible family of domains of $\VV'$. 
  Let $D'$ be an innermost island associated to the family $\II'$, and
let  $\JJ'=\II'|D$, $\phi=\phi_{D'}$.
For $j\in \JJ'$, let  us define $i(j)$ by the property
 $\phi(V_j')\subset U_{i(j)}$, and let
 $\II=\{i(j) : j\in \JJ'\}\cup\{0\}$.
 Then $\{U_i\}_{i\in \II}$ is an admissible family of domains of $\UU$, and
$$
  \sigma(\II'|D')\geq 
           {1\over 2}\left( (|\JJ'|-s)\nu +  s\; \mod P_0+
          \sum_{j\in \JJ',\, i(j)\not=0}\mod P_{i(j)}\right), 
$$
where $s=\#\{j: i(j)=0\}$ is the number of $\phi$-precritical pieces,
and $P_i$ are the maximal annuli enclosing $U_i$ in $\Delta$ rel $\II$.
\end{lem}

Using the transit map $\Psi$ from \S \ref{gpr}, we see that
the annuli $P_{i(j)}$ from Lemma \ref{main} 
are univalent pull-backs of some off-critical annuli
$R_{k(j)}$ of the family $\VV$ in $\De$. 
Incorporating it into  Lemma \ref{main}, we obtain:
$$ \si(I')\geq \si(I'| D')\geq  ( {1\over 2}\nu-\mu) +  \mu + 
       {1\over 2} \sum \mod R_{k(j)} \geq $$
\begin{equation}\label{definite}
\geq    ( {1\over 2}\nu-\mu) +  (\mu-\mod(R_0)) + \bar\si(\VV)  \geq  \bar\si(\VV)
\end{equation}
%
  Taking the infimum over all admissible families $I'$, we see that 
\begin{equation}\label{si} 
  \bar\si(\VV')\geq \bar\si(\VV).
\end{equation}

Consider now the principal nest $g_n: \cup V_i^n\ra \De^n$ of generalized quadratic-like maps.
The next map of the nest is related to the previous one  by the usual or parabolic generalized 
renormalization. We assume that if $g_{n+1}$ is the parabolic renormalization of $g_n$,
then the previous renormalization level  is non-central, i.e., 
 $g_{n-1}(0)\not\in  V^{n-1}_0$ (for otherwise $g_n|V^n_0 = g_{n-1} |V^n_0$
 and we could skip $g_n$ from  the  nest).
 Denote by $\mu_n= \mod(\De_n\sm V_n^0)$ the principal moduli of the $g_n$
 and by $\bar \si_n$ their asymmetric moduli.

\begin{cor}
  The asymmetric moduli $\bar\si_n$ stay away from 0.
  The principal moduli $\mu_n= \mod(\De_n\sm V_n^0)$ stay away from 0 for all parabolic maps
   $g_n$ .
\end{cor}

\begin{pf} The first statement immediately follows from (\ref{si}) and the corresponding 
 statement for the usual renormalization (Corollary 5.5 of \cite{puzzle}).

 Take a parabolic map $g_n$. If it is obtained from $g_{n-1}$  by the usual
generalized renormalization,
then the previous renormalization level  is non-central. 
  By Corollary 5.3 of \cite{puzzle},
$\mu_n\geq {1\over 2} \bar\si_{n-1}$. Hence these principal moduli $\mu_n$ stay away from 0.    

If $g_n$ is obtained from $g_{n-1}$ by the parabolic renormalization then 
by Lemma \ref{mu}, $\mu_n\geq \mu_{n-1}.$ Putting these together, we obtain the assertion.
\end{pf}



Let us go back to the estimate (\ref{definite}).
 If the hyperbolic distance between the domains $V_0$ and  $V_{k(j)}$  
(in the hyperbolic metric of $\De$) is bounded by some $d $, 
then $\mu\geq \mod(R_0) +\alpha(d)$, where $\alpha(d)>0$, 
and we observe a definite increase of the asymmetric modulus: 
  $$\si(\II')\geq \bar\si(\VV)+\alpha(d).$$ 

Since the tower under consideration has essentially bounded geometry, 
 the intervals $I_k$ stay a bounded hyperbolic distance from $I_0$ in the
hyperbolic metric of $J$. If at the same time the hyperbolic distance from some
$V_{k(j)}$ to $V_0$ in $\De$ is big ($\geq d$), then  
 the curve $\di \De$ must be pinched near the real line. In this case
we have a big adding in the  Gr\"otzsch inequality when we estimate $\nu$ by $2\mu$
(see \S\S 13.2- 13.3 in  \cite[Appendix A]{puzzle}): $\nu\geq 2\mu + \beta(d)$,
where $\beta(d)\to \infty$ as $d\to \infty$.
This  yields  a big increase of the asymmetric modulus:
   $$\si(\II')\geq \bar\si(\VV)+{1\over 2} \beta(d).$$

\msk
Thus the asymmetric moduli $\bar\si_n$ grow. Since it is incompatible with 
the essentially bounded geometry on the real line, we arrive at a contradiction.


\begin{thebibliography}{*****}

\bibitem[D1]{D1} A. Douady. Does a Julia set depend continuously on the polynomial?
   In: ``Complex dynamical systems'', Proceed. Symp. Apllied Math., v. 49, AMS,
   91 - 138. 

\bibitem[D2]{D2} A. Douady (with the participation of X.~Buff, R.~Devaney
  and P.~Sentenac). Baby Mandelbrot sets are born in cauliflowers. 
  In the book ``The Mandelbrot set, themes and variations'',
  London Math. Soc. Lect. Note Series, v. 274 (2000), 19-36.    

\bibitem[DH]{DH}
 A. Douady \& J.H. Hubbard. On the dynamics of polynomial-like maps.
    Ann. Sc. \'{E}c. Norm. Sup., v. 18 (1985), 287-343.
   
\bibitem[Hi]{Hi} B. Hinkle. Parabolic limits of renormalization.
  Preprint IMS at Stony Brook \# 1997/7. To appear in Erg. Th. \& Dynam. Syst.

\bibitem[L1]{attrac} M. Lyubich. Combinatorics, geometry and attractors of
  quasi-quadratic maps. Ann. Math., v. 140 (1994), 347 - 404. 

\bibitem[L2]{puzzle} M. Lyubich. Dynamics of quadratic polynomials, I-II. 
   Acta Math., v. 178 (1997), 185-297. 

\bibitem[L2]{parapuzzle} M. Lyubich. Dynamics of quadratic polynomials, III.
  Parapuzzle and SBR measure. 
Asterisque volume in honor of Douady's 60th birthday ``G\'eom\'etrie complexe et syst\'emes
 dynamiques'', v. 261 (2000),  173 - 200. 

\bibitem[LY]{LY} M. Lyubich \& M. Yampolsky. Dynamics of quadratic polynomials:
Complex bounds for real maps. Ann. Institut Fourier, v. 47 (1997), 1219 - 1254. 

\bibitem[M]{M} J. Milnor. On the concept of attractor. Comm. Math. Physics.,
    v. 99 (1985), 177-195, and v. 102 (1985), 517-519. 

\bibitem[McM]{McM} C. McMullen. Complex dynamics and renormalization. 
   Annals of Math. Studies, \# 135 (1994), Princeton Univ. Press. 

\bibitem[MS]{MvS} W. de Melo \& S. van Strien. One dimensional dynamics.
     Springer Verlag, 1993. 

\bibitem[Sh]{Sh} M. Shishikura. The parabolic bifurcation of rational maps.
  Col\'{o}qio 19 Brasileiro  de Matem\'{a}tica.
\end{thebibliography}
\end{document}